\documentclass[12pt, reqno]{amsart}
\usepackage{amsmath}
\usepackage{amsfonts}
\usepackage{amssymb}
\usepackage{amsthm}
\usepackage{mathtools}
\usepackage{mathrsfs}
      \theoremstyle{definition}
\newtheorem{defi}{Definition}[section]
\theoremstyle{plain}
\newtheorem{thm}[defi]{Theorem}

\newtheorem{cor}[defi]{Corollary}
\newtheorem{lemma}[defi]{Lemma}
\theoremstyle{definition}
\newtheorem{rmk}[defi]{Remark}
\newtheorem{ex}[defi]{Example}

\newtheorem*{pf}{Proof}
\newtheorem*{ack}{Acknowledgement}
\theoremstyle{plain}
\newtheorem{fact}[defi]{Fact}
\newtheorem{claim}[defi]{Claim}

\newcommand{\Z}{\mathbb{Z}}
\newcommand{\Q}{\mathbb{Q}}
\newcommand{\R}{\mathbb{R}}

\newcommand{\Div}{\text{Div}}

\newcommand{\supp}{\text{supp}}

      \makeatletter
      \def\@setcopyright{}
      \def\serieslogo@{}
      \newenvironment{proofof}[1]{\par
  \pushQED{\qed}%
  \normalfont \topsep6\p@\@plus6\p@\relax
  \trivlist
  \item[\hskip\labelsep
        \bfseries
    Proof of #1\@addpunct{.}]\ignorespaces
}{%
  \endtrivlist\@endpefalse
}
      \makeatother
   
\begin{document}
   \author{Roberto Laface}
   \address{Department of Mathematics "Federigo Enriques", University of Milan, Via Saldini 50, 20133 Italy}
   \email{roberto.laface@studenti.unimi.it}
   \title[Zariski Decompositions]{On Zariski decomposition\\with and without support}
   \begin{abstract}
     In this paper we study Zariski Decomposition with support in a negative definite cycle, a variation introduced by Y.\,Miyaoka. We provide two extensions of the original statement, which was originally meant for effective $\Q$-divisors: we can either state it for any $\Q$-divisor, or we can take the support to be in any cycle. Ultimately, we present a new approach to Zariski Decomposition of pseudo-effective $\Q$-divisors, which consists in iterating Zariski Decomposition with support.
   \end{abstract}
   \maketitle
	\thispagestyle{empty}
	\section*{Introduction}
		In 1962, O.\,Zariski introduced Zariski Decomposition for effective divisors in his paper \cite{zar}: given an effective divisor $D$ on a nonsingular projective surface $X$, we can write
		\[D=P+N,\]
		where $P$ and $N$ are effective, $P$ is nef, $N$ has negative definite intersection matrix, and $P \perp N$ with respect to the intersection product. In 1979, T.\,Fujita \cite{fuji} extended Zariski's result to the slightly larger cone of pseudo-effective divisors: the same result holds for pseudo-effective divisors, but the nef part $P$ is not necessarily effective anymore. Thirty years later, Y.\,Miyaoka \cite{miya} introduced the concept of Zariski Decomposition with support in a negative definite cycle, in a context which was far from where the actual problem was born: given a negative definite cycle $G$, any effective divisor $D$ can be written as 
		\[D= P_G+N_G,\]
		with the same properties as before, with an exception made for $P_G$, which now is $G$-nef only and is in general not effective. Meanwhile, in the same year, T.\,Bauer provided a simpler proof of Zariski's result, which was based on a maximality argument relative to the nef part of the given divisor, rather than on the sophisticated procedure Zariski used to build the negative one.\\
			
The aim of this paper is to connect the works we have presented in this timeline. First of all, we recall Zariski Decomposition with support in a negative definite cycle in the sense of Miyaoka (\cite{miya}, Proposition 2.1), and we use Bauer's method \cite{bauer} to give a new proof of this result, which Miyaoka states only, referring instead to Zariski's argument in \cite{zar} (Theorem \ref{gzardec}). Then, we doubly generalize Miyaoka's result: on one hand, we extend it to any divisor, once again using Bauer's idea (Theorem \ref{gzardecps}), and on the other, we prove the result for pseudo-effective divisors with the relaxed hypothesis that the support is in any cycle, not necessarily negative definite (Theorem \ref{gzardecnew}). These conditions cannot be dropped at once, as an elementary counterexample shows (Remark \ref{counterex}). Finally, we will provide a new proof of Zariski Decomposition of pseudo-effective divisors (\cite{fuji}, Theorem 1.12), which is obtained as an application of Zariski Decomposition with support in a negative definite cycle (Theorem \ref{zafudec}).
	
	\section{A new proof of Zariski Decomposition\\with support in a negative definite cycle}\label{sec1}
	Let $X$ be a surface, i.e.\,a 2-dimensional nonsingular projective variety over an algebraically closed field $k$.
\begin{defi}
Given a divisor $D =\sum_{i=1}^n d_i D_i$, where $d_i \in \Z$ and $D_i$ is an irreducible curve, $\forall i=1, \dots, n$, the matrix
\[\mu_D:=\begin{bmatrix}
D_1 . D_1 & \cdots & D_1 . D_q \\ \vdots & \ddots & \vdots \\ D_q . D_1 & \cdots & D_q . D_q
\end{bmatrix}\]
is called the \emph{intersection matrix of} $D$.
\end{defi}
\noindent As a consequence, the intersection matrix of a divisor is independent on the coefficients of the irreducible components. Therefore, the definition above can be obviously extended to divisors with rational (or real) coefficients; divisors with rational coefficients are called $\Q$-divisors, and we denote their vector space by $\Q\Div(X)$.\newline
Notice also that every $\Q$-divisor $D=\sum_{i=1}^n d_i D_i$ induces a quadratic form $\Phi_D$ on $\Q^n$ through its intersection matrix.
\begin{align*}
\Phi_D \colon &\Q^n \longrightarrow \Q \\
&\underline{v} \longmapsto \underline{v}^T \mu_D \underline{v}
\end{align*}
\begin{defi}
A finite sum $G=\sum_{i=1}^m G_i$ of irreducible curves $G_i \subset X$ is said to be a \textit{negative definite cycle} if the intersection matrix $\mu_G$ is negative definite, meaning that the quadratic form $\Phi_G$ induced by $G$ is negative definite.
\end{defi}

\begin{rmk}
From the definition, it follows that the components $G_i$ of a negative definite cycle must be distinct, i.e.\,$G$ is reduced: if not, the matrix $\mu_G$ would have 2 equal columns, hence the quadratic form $\Phi_G$ would not be negative definite.
\end{rmk}

\begin{ex}
The $f$-exceptional locus of a surjective morphism of surfaces $f: X \longrightarrow Y$ (the union of the curves that $f$ contracts to a point) is a typical example of negative definite cycle.
\end{ex}

\begin{ex}
Consider any $\Q$-divisor $D \subset X$ and write its Zariski Decomposition $D=P+N$, $N=\sum_i \nu_i N_i$; then $N_{\text{red}}:=\sum_i N_i$ is a negative definite cycle.
\end{ex}

\begin{rmk}\label{gsupp}
Let $D \in \Q \Div (X)$ be a divisor which is supported on a negative definite cycle $G=\sum_{i=1}^m G_i$, i.e. $D=\sum_{i=1}^m d_i G_i$, $d_i \in \Q \ \forall i=1, \dots , m$. Then
\[D^2=\Big( \sum_{i=1}^m d_i G_i \Big)^2 = \sum_{i,j=1}^m d_i d_j (G_i . G_j)=
\begin{pmatrix}
d_1 & \dots & d_m
\end{pmatrix}
\mu_G
\begin{pmatrix}
d_1\\
\vdots\\
d_m
\end{pmatrix} \leq 0,\]
and equality holds if and only if $D=0$. 
\end{rmk}

\begin{defi}
A $\mathbb{Q}$-divisor $D$ is \textit{$G$-nef} (i.e.\,numerically effective on $G$) if $D . G_i \geq 0$ for every $i$.
\end{defi}

Now, we introduce two new orderings:

\begin{enumerate}
\item given $C=\sum_{i=1}^n c_i C_i, \, D=\sum_{i=1}^n d_i C_i \in \Q\Div (X)$, we set $C \leq D$ if and only if $c_i \leq d_i \ \forall i=1, \dots ,n$;
\item given $\underline{x}=(x_1, \dots , x_n), \ \underline{y}=(y_1, \dots , y_n) \in \Q^n$, we put $\underline{x} \leq \underline{y}$ if and only if $x_i \leq y_i \ \forall i =1, \dots , n$.
\end{enumerate}
We will use these two orderings in the proof of the next result.

\begin{thm}[Zariski Decomposition with support in a negative definite cycle, Proposition 2.1 of \cite{miya}]\label{gzardec} Let $G=\sum_{i=1}^m G_i \subset X$ be a negative definite cycle and let $D$ be an effective $\mathbb{Q}$-divisor on $X$. Then there exists a unique decomposition $D=P+N$ into $\mathbb{Q}$-divisors which satisfies the following conditions:

\begin{enumerate}

\item[(a)] both $P$ and $N$ are effective;

\item[(b)] $N$ is supported on a subset of $G$, i.e. $N=\sum_i \nu_i G_i$, $\nu_i \geq 0$;

\item[(c)] $P$ is $G$-nef;

\item[(d)] $P$ and $N$ are mutually orthogonal, i.e. $P . N = 0$ (hence $D^2 = P^2 + N^2$ and, in view of \text{\textnormal{(c)}}, $P$ is numerically trivial on $N$, i.e. $P . G_i =0$ for each $G_i \subset supp \, N$). 

\end{enumerate}

\noindent Furthermore, $P$ is the largest effective $\mathbb{Q}$-subdivisor of $D$ that is $G$-nef:

\begin{enumerate}

\item[(e)] if a $\mathbb{Q}$-divisor $E$ with $0 \leq E \leq D$ is $G$-nef, then $E \leq P$.

\end{enumerate}

\end{thm}

The so-obtained decomposition is called $G$-decomposition of $D$, $P$ is the $G$-nef part of $D$, while $N$ is said to be the $G$-negative part. Miyaoka's paper \cite{miya} is lacking in any proof of this result, and rather it refers to \cite{zar} for the proof. Instead, we will give a proof which provides a concrete application of Bauer's proof in \cite{bauer}, distancing itself from the original idea of Zariski.

\begin{pf}
We start by proving the existence of such a decomposition. Let $D=\sum_{i=1}^n d_i D_i$, with $D_i$ irreducible curve and $d_i \in \mathbb{Q}_{> 0}$, $\forall i=1, \dots ,n$. Consider now $A$ such that $0 \leq A \leq D$, $A = \sum_{i=1}^n x_i D_i$, $0 \leq x_i \leq d_i$. We now have that
\begin{align*}
\text{$A$ is $G$-nef \ $\Longleftrightarrow$} & \ \text{$A \cdot G_j \geq 0 \ \forall j=1, \dots, m$} \\
\text{$\Longleftrightarrow$} & \ \text{$\sum_{i=1}^n x_i D_i \cdot G_j \geq 0 \ \forall j=1, \dots, m$.} \tag{*}
\end{align*}
\begin{claim}\label{cl1} This system of inequalities for the rational numbers $x_i$ has a maximal solution (with respect to the ordering $\leq$) in the rational cuboid $[0, d_1] \times \dots \times [0, d_n] \subset \mathbb{Q}^n$. \end{claim}
\begin{pf}
Indeed, the subset $K$ of $\mathbb{Q}^n$ described by these inequalities is the intersection of finitely many half-spaces; notice that we always have a solution (the vector $\underline{x}=\underline{0}$). Consider the family of hyperplanes \[H_t := \big\lbrace (x_1 , \dots , x_n) \in \mathbb{Q}^n \colon \sum_{i=1}^n x_i = t \sum_{i=1}^n d_i \big\rbrace;\] then, there is a maximal $t$ such that $H_t$ intersects $K$ (and the point of intersection is a vertex of $K$). \qed
\end{pf}
\noindent Now let $P$ be a $\Q$-divisor defined by a maximal solution to the system of inequalities above, $P=\sum_{i=1}^n b_i D_i$. Set $N:=D-P$; then conditions ($a$) and ($c$) are satisfied by construction. We now prove ($b$) and ($d$): notice that in case $D=P$, the last two conditions hold, and thus we can assume that $N$ is nonzero.
\begin{enumerate}
\item[$(b)$] Write $N = \sum_i \nu_i D_i$. For a fixed $i$, consider the intersections $D_i . G_j$, $j=1, \dots , m$; if $D_i . G_j \geq 0$ $\forall j=1, \dots , m$, then $P + \varepsilon D_i$ is $G$-nef (for a suitable $\varepsilon > 0$), and this contradicts the maximality of $P$. Thus there must exist a $j=j(i)$ such that $D_i = G_j$ and $D_i ^2 = G_j ^2 <0$. Since this holds for every $i$, we get $N = \sum_i \nu_i G_i$ after rearranging indexes. 
\item[$(d)$] By $(b)$, $N = \sum_i \nu_i G_i$. If $P . G_i >0$, with $G_i \subseteq \supp (N)$, then $P + \varepsilon G_i \leq D$ and $P + \varepsilon G_i$ is $G$-nef, for small enough $\varepsilon >0$, contradicting the maximality of $P$. Then $P.G_i=0$, because $P$ is $G$-nef, and $P.N=0$.
\end{enumerate}
Assume now that we are given a decomposition $D=P+N$.
\begin{claim}\label{cl2}
$P$ is a maximal $G$-nef subdivisors of $D$.
\end{claim}
\begin{pf}
Indeed, given a $G$-nef divisor $P': \ P \leq P' \leq D=P+N$, we have that $P'=P+ \sum_{i \in \mathfrak{I}} y_i G_i$, being
\[\mathfrak{I}:= \lbrace k : \ G_k \subseteq \supp (N) \rbrace.\] 
$G$-nefness of $P'$ and orthogonalilty between $P$ and $N$ imply that
\[0 \leq P'. G_j = \sum_{i \in \mathfrak{I}} y_i (G_i.G_j) \quad \forall j \in \mathfrak{I}.\]
Then, by multiplying by $y_j$, we get $\sum_{i \in \mathfrak{I}} y_iy_j (G_i.G_j) \geq 0$ $\forall j \in \mathfrak{I}$, and by summing over all $j \in \mathfrak{I}$, we get
\begin{align*}
0 &\leq \sum_{i,j \in \mathfrak{I}} y_iy_j(G_i.G_j)= \sum_{i \in \mathfrak{I}} y_i^2 G_i^2 + \sum_{i\neq j} y_iy_j(G_i.G_j) \leq \\
& \leq \sum_{i \in \mathfrak{I}} y_i^2 G_i^2 + 2 \sum_{i\neq j} y_iy_j(G_i.G_j) = \Big( \sum_{i \in \mathfrak{I}} y_iG_i \Big)^2= \Phi_G (\underline{y}),
\end{align*}
where $\underline{y} \in \Q^m$ is the vector whose components are the coefficients of the $G_i$'s ($y_i=0 \ \forall i: \, i \notin \mathfrak{I}$). Since $\Phi_G$ is negative definite, it can only be $\Phi_G(\underline{y})=0$, and this happens if and only if $\underline{y}=\underline{0}$, yielding $P'=P$ and thus maximality of $P$.
\qed
\end{pf}
\noindent Now we are left to prove uniqueness. We show that a maximal $G$-nef $\mathbb{Q}$-subdivisor of $D$ is in fact unique.
\begin{lemma}\label{lemma1}
If $P' = \sum_{i=1}^n x'_i D_i$ and $P'' = \sum_{i=1}^n x''_i D_i$ are $G$-nef $\mathbb{Q}$-subdivisors of $D$, then so is $P= \max(P',P''):=\sum_{i=1}^n x_i D_i$, where $x_i := max(x'_i, x''_i)$. 
\end{lemma}
\begin{pf}
Showing that $P$ is $G$-nef is equivalent to showing that it is $G_i$-nef for every $i=1, \dots , m$. If $G_i \nsubseteq \supp (P)$, then
\begin{align*}
(P-P'). G_i &= \sum_{k=1}^n \lbrace \max (x'_k,x''_k)-x'_k \rbrace \underbrace{(D_k.G_i)}_{\geq 0} \geq 0,
\end{align*}
and thus $P . G_i = (P-P'). G_i + P'. G_i \geq 0$; otherwise, $G_i =D_j$ for some $j$. Without loss of generality we can assume that $x'_j \geq x''_j$; now we get
\begin{align*}
(P-P'). G_i = (P-P').D_j &= \sum_{k=1}^n \lbrace \max (x'_k,x''_k)-x'_k \rbrace (D_k.D_j) =\\
&=\sum_{k \neq j} \max (x'_k,x''_k)-x'_k \rbrace (D_k.D_j) \geq 0,
\end{align*}
and thus 
\[P.G_i=\underbrace{(P-P').G_i}_{\geq 0}+ \underbrace{P'.G_i}_{\geq 0} \geq 0,\]
i.e. $P$ is $G_i$-nef. Since this holds for every $i$, we are done.
\qed
\end{pf}
\noindent By the lemma, $P$ is the maximal $G$-nef subdivisor of $D$; hence ($e$) is proven, and this concludes the proof of the theorem.
\qed
\end{pf}

\begin{rmk}
Assume we are given a Zariski Decomposition $D=P+N$ with support in a negative cycle $G$; then the negative part $N=\sum_i \nu_i G_i$ has negative definite matrix because $\Phi_N$ is the restriction of $\Phi_G$ to the subspace $V$ of $\Q^m$ defined by
\[V:=\lbrace \underline{x} \in \Q^m \vert x_i=0 \ \forall i: \,\nu_i=0 \rbrace.\]
\end{rmk}
	
	\section{Two improvements of Miyaoka's result}
	
In this section we provide two generalizations of Theorem \ref{gzardec}. We first extend the original result to any $\Q$-divisor on a surface.
\begin{thm}[Zariski Decomposition with support in negative definite cycle for $\Q$-divisors]\label{gzardecps}
Let $G=\sum_{i=1}^q G_i$ be a negative definite cycle and let $D$ be a $\mathbb{Q}$-divisor on $X$. Then there exists a unique decomposition $D=P+N$ into $\mathbb{Q}$-divisors which satisfies the following conditions:
\begin{enumerate}
\item[(a)] $N$ is effective;
\item[(b)] N is supported on a subset of $G$, i.e. $N=\sum \nu_i G_i$, $\nu_i \geq 0$;
\item[(c)] $P$ is $G$-nef;
\item[(d)] $P$ and $N$ are mutually orthogonal, i.e. $P . N = 0$ (hence $D^2 = P^2 + N^2$ and, in view of (c), $P$ is numerically trivial on $N$, i.e. $P . G_i =0$ for each $G_i \subset \supp \, N$). 
\end{enumerate}
Furthermore, $P$ is the largest effective $\mathbb{Q}$-subdivisor of $D$ that is nef on $G$:
\begin{enumerate}
\item[(e)] if a $\mathbb{Q}$-divisor $E$ with $0 \leq E \leq D$ is nef on $G$, then $E \leq P$.
\end{enumerate}
\end{thm}
The key lemma for proving the theorem is the following
\begin{lemma}
Let $D$ be a $\Q$-divisor, and let $G=\sum_{i=1}^q G_i$ be a negative definite cycle. Then there exists a subdivisor $P \leq D$ such that $P$ is $G$-nef.
\end{lemma}
\begin{pf}
If $D$ is $G$-nef, we set $P:=D$; otherwise, $D$ is negative on some of the $G_i$'s. Define
\[P \equiv P(\underline{x}):= D - \sum_{i=1}^q x_i G_i,\]
where $\underline{x}=(x_1 , \dots , x_q) \in \Q_{\geq 0}^q$. Then $P$ is $G$-nef if and only if
\begin{gather*}
P . G_j \geq 0 \quad \forall j=1, \dots , q \ \Longleftrightarrow \sum_{i=1}^q x_i (G_i . G_j) \leq D . G_j \quad \forall j=1, \dots , q 
\end{gather*}
and the last condition is equivalent to the matrix inequality
\begin{align}\label{matrix1}
\begin{bmatrix}
G_1 . G_1 & \cdots & G_1 . G_q \\ \vdots & \ddots & \vdots \\ G_q . G_1 & \cdots & G_q . G_q
\end{bmatrix}
\begin{bmatrix}
x_1 \\ \vdots \\ x_q 
\end{bmatrix} \leq
\begin{bmatrix}
D . G_1 \\ \vdots \\ D . G_q
\end{bmatrix}.
\end{align}
The inequality \eqref{matrix1} is equivalent to the following homogeneous one
\begin{align}\label{matrix2}
\begin{bmatrix}
-G_1 . G_1 & \cdots & -G_1 . G_q &  D . G_1  \\ \vdots & \ddots & \vdots  & \vdots \\ -G_q . G_1 & \cdots & -G_q . G_q & D . G_q \\ 0 & \cdots & 0 & 1
\end{bmatrix}
\begin{bmatrix}
x_1 \\ \vdots \\ x_q \\ 1
\end{bmatrix}
\geq 0;
\end{align}
however \eqref{matrix2} has a solution if and only if there exists a solution to the system
\begin{align}\label{matrix3}
\begin{bmatrix}
-G_1 . G_1 & \cdots & -G_1 . G_q &  D . G_1  \\ \vdots & \ddots & \vdots  & \vdots \\ -G_q . G_1 & \cdots & -G_q . G_q & D . G_q \\ 0 & \cdots & 0 & 1
\end{bmatrix}
\begin{bmatrix}
x_1 \\ \vdots \\ x_q \\ t
\end{bmatrix}
\geq 0;
\end{align}
In fact, it sufficies to divide or multiply by $t$ a solution to one of the two systems in order to get a solution to the other one; notice that if we denote the matrix in \eqref{matrix3} by $M$, then $M$ is such that all principal minors of $M$ are nonnegative: for, notice that $M$ is built out of $- \mu_G$, which is positive definite since $\mu_G$ is negative definite. The claim now follows applying Laplace's theorem for computing determinants to the last row.\newline
Finally, we get the result by applying
\begin{fact}[Theorem 1 and Remark in \cite{kara}]
Let $A$ be an $n \times n$ real matrix.
\begin{enumerate}
\item[(i)] If all the principal minors are positive, then the system
\begin{align}\label{syst1}
\begin{cases}
\underline{x}\geq 0 \\
A \underline{x} > 0
\end{cases}
\end{align}
has a solution.
\item[(ii)] If all the principal minors are nonnegative, then the system
\begin{align}
\begin{cases}
\underline{x}\geq 0 \\
A \underline{x} \geq 0
\end{cases}
\end{align}
has a solution.
\end{enumerate} 
Moreover, \eqref{syst1} has a solution if and only if the system
\begin{align}
\begin{cases}
\underline{x}> 0 \\
A \underline{x} > 0
\end{cases}
\end{align}
has a solution.
\end{fact}
\qed
\end{pf}

\begin{proofof}{Theorem \ref{gzardecps}}
The lemma ensures the existence of a $G$-nef subdivisor of $D$; now we can choose the solution $\underline{x}$ to be minimal with respect to the ordering $\leq$ defined in Section \ref{sec1}. This leads to a maximal subdivisor $P$ of $D$, with respect to the property of being $G$-nef. Thus, we can use the same maximality argument of Theorem \ref{gzardec} to get existence and uniqueness of the Zariski Decomposition with support in $G$.
\qed
\end{proofof}

A second extension of Miyaoka's Theorem \ref{gzardec} holds for pseudo-effective $\Q$-divisors: this variation allows us to take support in any cycle, not necessarily negative definite. The main idea is the following: given a pseudo-effective $\Q$-divisor, we fix its non-$G$-nefness step by step by iterating Theorem \ref{gzardecps}.\newline
We first recall some result we need to deal with pseudo-effective $\Q$-divisors.

\begin{lemma}[Lemma 1.8 of \cite{fuji}]\label{fuj9}
Let $\lbrace C_i \rbrace_{i=1, \dots , q}$ be a family of integral curves, and let $E:=\sum_{i=1}^q a_i C_i$ be a $\Q$-divisor. If $D$ is a pseudo-effective $\Q$-divisor such that $(D-E).C_i \leq 0$ $\forall i=1, \dots, q$, then $D-E$ is pseudo-effective.
\end{lemma}

\begin{lemma}[Lemma 1.10 of \cite{fuji}]\label{fuj7}
Let $\lbrace C_i \rbrace_{i=1, \dots , q}$ be a family of integral curves such that the matrix
\begin{align*}
\begin{bmatrix}
C_1 . C_1 & \cdots & C_1 . C_r \\ \vdots & \ddots & \vdots \\ C_r . C_1 & \cdots & C_r . C_r
\end{bmatrix}
\end{align*}
is negative definite for some $r < q$. If $D$ is a pseudo-effective $\Q$-divisor such that $D . C_i \leq 0$ for every $i=1, \dots , q$ and $D . C_i <0$ for $i=r+1, \dots , q$, then the matrix
\begin{align*}
\begin{bmatrix}
C_1 . C_1 & \cdots & C_1 . C_{r} & C_1 . C_{r+1} & \cdots & C_1 . C_{q}\\
\vdots & \ddots & \vdots & \vdots & \ddots & \vdots\\
C_{r} . C_1 & \cdots & C_{r} . C_{r} & C_{r} . C_{r+1} & \cdots & C_{r} . C_{q}\\
C_{{r}+1} . C_1 & \cdots & C_{{r}+1} . C_{{r}} & C_{{r}+1} . C_{r+1} & \cdots & C_{r+1} . C_{q}\\
\vdots & \ddots & \vdots & \vdots & \ddots & \vdots\\
C_{q} . C_1 & \cdots & C_{q} . C_{r} & C_{q} . C_{r+1} & \cdots & C_{q} . C_{q}
\end{bmatrix}
\end{align*}
is also negative definite.
\end{lemma}

We finally state and prove the following
\begin{thm}[Zariski Decomposition with support in a cycle for pseudo-effective $\Q$-divisors]\label{gzardecnew}
Let $G=\sum_{i=1}^m G_i \subset X$ be a cycle and let $D$ be a pseudo-effective $\mathbb{Q}$-divisor on $X$. Then there exists a unique decomposition $D=P+N$ into $\mathbb{Q}$-divisors which satisfies the following conditions:

\begin{enumerate}

\item[(a)] $N$ is effective;

\item[(b)] $N$ is supported on a subset of $G$, i.e. $N=\sum_i \nu_i G_i$, $\nu_i \geq 0$;

\item[(c)] $P$ is $G$-nef;

\item[(d)] $P$ and $N$ are mutually orthogonal, i.e. $P . N = 0$ (hence $D^2 = P^2 + N^2$ and, in view of \text{\textnormal{(c)}}, $P$ is numerically trivial on $N$, i.e. $P . G_i =0$ for each $G_i \subset supp \, N$). 

\end{enumerate}

\noindent Furthermore, $P$ is the largest effective $\mathbb{Q}$-subdivisor of $D$ that is $G$-nef:

\begin{enumerate}

\item[(e)] if a $\mathbb{Q}$-divisor $E$ with $0 \leq E \leq D$ is $G$-nef, then $E \leq P$.

\end{enumerate}

\end{thm}

\begin{pf}
We proceed by iterating Theorem \ref{gzardecps}. If $D$ is $G$-nef, then we are done; otherwise, define
\[G^{(1)}:= \sum_{i=1}^{q_1} G_i,\]
where $D.G_i < 0$ $\forall i=1, \dots, q_1$; by Lemma \ref{fuj7}, $G^{(1)}$ is a negative definite cycle (apply to the case $r=0$). Hence, we can write the Zariski Decomposition with support in $G^{(1)}$ of $D$:
\[D=P_{G^{(1)}}+N_{G^{(1)}}.\]
If $P_{G^{(1)}}$ is $G$-nef, then we are done; otherwise, we notice that $P_{G^{(1)}}$ is pseudo-effective by Lemma \ref{fuj9}. Now, define
\[G^{(2)}:= \sum_{i=1}^{q_2} G_i,\]
where $P_{G^{(1)}}.G_i < 0$ $\forall i=q_1+1, \dots, q_2$ (and $P_{G^{(1)}}.G_i=0$ $\forall i=1, \dots, q_1$); again by Lemma \ref{fuj7}, $G^{(2)}$ is a negative definite cycle. Write the Zariski Decomposition with support in $G^{(2)}$
\[P_{G^{(1)}}=P_{G^{(2)}}+N_{G^{(2)}},\]
and obtain the decomposition
\[D=P_{G^{(2)}}+ \big( N_{G^{(1)}}+N_{G^{(2)}} \big).\]
If $P_{G^{(2)}}$ is $G$-nef, then we are done; otherwise we repeat this process, which must come to an end since $G$ is a finite sum of integral curves.
\qed
\end{pf}

\begin{rmk}\label{counterex}
Referring to Theorem \ref{gzardec}, we cannot relax both the hypothesis on the given divisor and the cycle at the same time.\\
There is an elementary counterexample in this regard: we can take $X=\mathbb{P}^2_{\mathbb{C}}$, $D=-H$ and $G=H$, $H$ being the hyperplane divisor. Then, $-H.H=-1$, and thus $-H$ is not $H$-nef, but also $H^2>0$, meaning that $H$ is not a negative definite cycle. Hence, there is no Zariski Decomposition with support in $H$ of $-H$.
\end{rmk}
\section{A new approach to\\Zariski Decomposition of pseudo-effective $\Q$-divisors}

In \cite{fuji}, Fujita generalized the idea of Zariski Decomposition to pseudo-effective divisors. In this section, we show an alternative approach to Fujita's result. Before moving to the main result, we recall the following

\begin{cor}[Corollary 1.11 of \cite{fuji}]\label{fuj8}
Let $D$ be a pseudo-effective $\Q$-divisor. Then there exists only a finite number of integral curves on which $D$ is negative, meaning that there are only finitely many integral curves $C$ such that $D . C<0$.
\end{cor}

\begin{thm}[Zariski Decomposition of pseudo-effective $\Q$-divisors]\label{zafudec} Let $D$ be a pseudo-effective $\Q$-divisor. Then there are uniquely determined $\mathbb{Q}$-divisors $P$ and $N$ (called nef and negative part respectively) with 

\begin{align*}
D=P+N
\end{align*}

such that:

\begin{enumerate}
\item[(i)] $P$ is nef;
\item[(ii)] $N$ is effective 
\item[(iii)] $N$ is either zero or it has negative definite intersection matrix;
\item[(iv)] $P . C =0$ for every irreducible component $C$ of $N$.
\end{enumerate}
\end{thm}

\begin{pf}
If $D$ is nef, then set $P:=D$ and $N:=0$. Otherwise, there are integral curves on which $D$ is negative, and these are in a finite number by Corollary \ref{fuj8}; let this family be $\mathscr{C}^{(1)}:=\lbrace C_1, \dots , {C_q}_1 \rbrace$. The matrix
\begin{align*}
\begin{bmatrix}
C_1 . C_1 & \cdots & C_1 . C_{q_1} \\ \vdots & \ddots & \vdots \\ C_{q_1} . C_1 & \cdots & C_{q_1} . C_{q_1}
\end{bmatrix}
\end{align*}
is negative definite by Lemma \ref{fuj7} (apply to the case $r=0$) and hence
\[G:=\sum_{i=1}^{q_1} C_i\]
is a negative definite cycle; consequently, we can take the Zariski Decomposition of $D$ with support in $G$
\[D=P_G+N_G.\]
Notice that $C_i$ appears in $N_G$ with nonzero coefficient, $\forall i=1, \dots, q_1$: indeed, if $C_i \nsubseteq \supp \, N_G$, then we get
\[0>D . C_i =(P_G+N_G) . C_i = \underbrace{P_G . C_i}_{\geq 0} + \underbrace{N_G . C_i}_{\geq 0} \geq 0,\]
which contradicts the negativity of $D$ on $C_i$; moreover, $P_G . C_i =0$ for every $i=1, \dots ,q_1$, because of the properties of the Zariski Decomposition with support. Now, if $P_G$ is nef, then we can put $P:=P_G$ and $N:=N_G$; otherwise, we notice that $P_G$ is pseudo-effective (Lemma \ref{fuj9}). Hence, there are curves $\mathscr{C}^{(2)}:=\lbrace C_{{q_1}+1} , \dots, C_{q_2} \rbrace$ on which $P_G$ is negative, and these are such that $\mathscr{C}^{(1)}\cap \mathscr{C}^{(2)}=\emptyset$ because $P_G$ is $G$-nef. Lemma \ref{fuj7} implies that the matrix
\begin{align*}
\begin{bmatrix}
C_1 . C_1 & \cdots & C_1 . C_{q_1} & C_1 . C_{q_1+1} & \cdots & C_1 . C_{q_2}\\
\vdots & \ddots & \vdots & \vdots & \ddots & \vdots\\
C_{q_1} . C_1 & \cdots & C_{q_1} . C_{q_1} & C_{q_1} . C_{q_1+1} & \cdots & C_{q_1} . C_{q_2}\\
C_{{q_1}+1} . C_1 & \cdots & C_{{q_1}+1} . C_{{q_1}} & C_{{q_1}+1} . C_{q_1+1} & \cdots & C_{q_1+1} . C_{q_2}\\
\vdots & \ddots & \vdots & \vdots & \ddots & \vdots\\
C_{q_2} . C_1 & \cdots & C_{q_2} . C_{q_1} & C_{q_2} . C_{q_1+1} & \cdots & C_{q_2} . C_{q_2}
\end{bmatrix}
\end{align*}
is negative definite and
\[G':=G+\sum_{i=q_1+1}^{q_2} C_i\]
is a negative definite cycle; now take the Zariski Decomposition of $P_G$ with support in $G'$:
\[P_G=P_{G'}+N_{G'}.\]
Once again, if $P_{G'}$ is nef, then we are done; otherwise, we continue applying this process, which has to stop after finitely many step because $\dim N^1(X)_\R<+\infty$. This concludes the proof of existence of the Zariski Decomposition.\newline
Assume we are given a Zariski Decomposition $D=P+N$ as in the statement, then
\begin{claim} $P$ is a maximal nef subdivisor of $D$.\end{claim}
\begin{pf}
For, assume that $P'$ is a nef divisor such that
\[P \leq P' \leq D=P+N;\]
then $P'=P+\sum_{k} \nu_k N_k$. By using nefness of $P'$ and orthogonality of $P$ and $N$, we get
\[0 \leq P' . N_j = \sum_k \nu_k (N_j . N_k) \ \forall j.\]
By multiplying by $\nu_j$, we have
\[0 \leq \sum_k \nu_j\nu_k (N_j . N_k) \ \forall j,\]
and summing over all $j$'s we get
\begin{align*}
0 &\leq \sum_{j,k} \nu_j \nu_k (N_j . N_k) = \sum_k \nu_k^2 N_k^2 + \sum_{j \neq k} \nu_j \nu_k (N_j . N_k) \leq \\
& \leq \sum_k \nu_k^2 N_k^2 + 2\sum_{j \neq k} \nu_j \nu_k (N_j . N_k) = (\sum_k \nu_k N_k)^2= \Phi_N(\underline{\nu}).
\end{align*}
However, $N$ has negative definite intersection matrix, hence $\Phi_N(\underline{\nu})=0$; thus $\underline{\nu}=\underline{0}$ and $P'=P$, i.e.\,$P$ is maximal.
\qed
\end{pf}
\noindent Uniqueness now follows by
\begin{lemma}\label{uniq}
If $P' = \sum_{i=1}^n x'_i D_i$ and $P'' = \sum_{i=1}^n x''_i D_i$ are nef $\mathbb{Q}$-subdivisors of $D$, then so is $P=\sum_{i=1}^n x_i D_i$, where $x_i := max(x'_i, x''_i)$. 
\end{lemma}
\begin{pf} $P$ is of course a $\Q$-subdivisor of $D$, so we only have to check nefness. Showing that $P$ is nef is equivalent to showing that it is $C$-nef for every integral $C$. If $C \nsubseteq \supp (P)$, then
\begin{align*}
(P-P'). C &= \sum_{k=1}^n \lbrace \max (x'_k,x''_k)-x'_k \rbrace \underbrace{(D_k.C)}_{\geq 0} \geq 0,
\end{align*}
and thus $P . C = (P-P'). C + P'. C \geq 0$; otherwise, $C =D_j$ for some $j$. Without loss of generality we can assume that $x'_j \geq x''_j$; now we get
\begin{align*}
(P-P'). C =(P-P').D_j &= \sum_{k=1}^n \lbrace \max (x'_k,x''_k)-x'_k \rbrace (D_k.D_j) =\\
&=\sum_{k \neq j} \max (x'_k,x''_k)-x'_k \rbrace (D_k.D_j) \geq 0,
\end{align*}
and thus 
\[P.C=\underbrace{(P-P').C}_{\geq 0}+ \underbrace{P'.C}_{\geq 0} \geq 0,\]
i.e. $P$ is $C$-nef. Since this holds for every $C$, we are done.
\qed
\end{pf}
\noindent The proof of the theorem is now complete.\qed
\end{pf}

\begin{rmk}
We notice that our procedure coincides with Fujita's at last: in fact, every divisor $N_G$ constructed in the proof via Zariski Decomposition with support is such that $N_G.G_i=D.G_i$, for all $G_i \subseteq \supp\,G$ ($P_G \perp N_G$), and this is exactly the condition Fujita used to build the divisors in his proof (see \cite{fuji}). Finally, there can only be one such divisor, since the matrix $\mu_G$ is negative definite and thus invertible. However, it is now clear how all the pieces are connected.
\end{rmk}
	
\begin{ack}
The present paper is part of the Master's Thesis of the author. It was written while the author was spending two academic semesters at the University of Bergen, funded by the Lifelong Learning Programme ERASMUS. The author is thankful to Prof. Andreas Leopold Knutsen (University of Bergen) and Prof. Antonio Lanteri (University of Milan), who supervised his work.
\end{ack}

\end{document}